\documentclass[12pt]{article}

\usepackage{amsmath}
\usepackage{amsfonts}
\usepackage{amsthm}
\usepackage{fancyvrb}

\newtheorem{theorem}{Theorem}

\newtheorem{property}{Property}
\newtheorem{definition}{Definition}
\newtheorem{problem}{Problem}

\numberwithin{equation}{section}
\numberwithin{figure}{section}

\title{A Golub-Welsch version for simultaneous Gaussian quadrature}
\author{Walter Van Assche\thanks{Supported by research project G0C9819N of FWO (Research Foundation -- Flanders).}\\ KU Leuven, Belgium}
\date{}  

\begin{document}
  \maketitle

\begin{abstract}
The zeros of type II multiple orthogonal polynomials can be used for quadrature formulas that approximate $r$ integrals of the same function $f$ with respect to $r$ measures $\mu_1,\ldots,\mu_r$
in the spirit of Gaussian quadrature. This was first suggested by Borges in 1994, even though he
does not mention multiple orthogonality. We give a method
to compute the quadrature nodes and the quadrature weights which extends the Golub-Welsch 
approach using the eigenvalues and left and right eigenvectors of a banded Hessenberg matrix.
This method was already described by Coussement and Van Assche in 2005 but it seems to
have gone unnoticed. We describe the result in detail for $r=2$ and give some examples.
\end{abstract}



\section{Introduction}
Suppose a real function $f$ is given and one wants to compute the integrals
\[   \int_\mathbb{R} f(x)\, d\mu_j(x)   , \qquad 1 \leq j \leq r  \]
for $r$ positive measures $\mu_1,\ldots,\mu_r$ with only $N$ function evaluations. In 1994 Borges \cite{Borges}
suggested Gauss-like quadrature rules and showed that one gets an optimal set of quadrature rules if one takes
the $N$ quadrature nodes at the zeros of a polynomial of degree $N$ which satisfies orthogonality conditions with
respect to the measures $\mu_1,\ldots,\mu_r$. Borges uses the terminology quasi-orthogonality but  such a polynomial is now known as a type II multiple orthogonal polynomial.
The quadrature rules are of the form
\[      \sum_{k=1}^N \lambda_{k,N}^{(j)} f(x_{k,N})  \approx  \int_\mathbb{R} f(x)\, d\mu_j(x),  \]
so that the integral for the measure $\mu_j$ corresponds to quadrature weights $\lambda_{k,N}^{(j)}$.
An important numerical aspect is then how to compute the quadrature nodes $\{x_{k,N}, 1 \leq k \leq N \}$
and the quadrature weights $\{ \lambda_{k,n}^{(j)}, 1 \leq k \leq N \}$ for $1 \leq j \leq r$ in an efficient and stable way. This problem was suggested in the survey on classical multiple orthogonal polynomials \cite[\S 4.3]{WVA-Couss}.
 The computation of the quadrature nodes and weights was first done by
Milovanovi\'c and Stani\'c for almost diagonal multi-indices \cite{Milovanovic} and later
this was extended for general multi-indices for the case of two or three quadrature rules by Tomovi\'c and Stani\'c \cite{Tomovic} and for arbitrary $r$ by Jovanovi\'c, Stani\'c and
Tomovi\'c \cite{Jovanovic}. The common nodes were calculated as eigenvalues of certain matrices, while the quadrature weights were calculated by solving certain systems of linear equations.
Simultaneous Gaussian quadrature for multiple orthogonal polynomials for two measures
$\mu_1$ and $\mu_2$ supported on two intervals $\Delta_1$ and $\Delta_2$ which are disjoint
(a so-called Angelesco system) was investigated in \cite{LubWVA}, and for multiple Hermite
polynomials (with $r=3$) in \cite{WVAVuer} but mainly from an analytical point of view, leaving out
the numerical aspects.

In 1969 Golub and Welsch \cite{GolubWelsch} showed that for Gaussian quadrature (the case when $r=1$)
one can use the Jacobi matrix for the corresponding orthogonal polynomials. The eigenvalues of the $N \times N$
Jacobi matrix are the quadrature nodes and the first components of the normalized eigenvectors are related to the
quadrature weights. This approach was extended to simultaneous Gaussian quadrature (in the sense of Borges) by
Coussement and Van Assche \cite{JCoussWVA} but this extension seems to have gone unnoticed. The role of the
Jacobi matrix is now played by a banded Hessenberg matrix for which the eigenvalues are the required quadrature nodes,
and one needs the right and the left eigenvectors to find the quadrature weights.
The goal of this paper is to give a simplified proof of this result for two measures $\mu_1$ and $\mu_2$. We will
give some examples involving measures $d\mu_1 = w_1(x)\,dx$ and $d\mu_2(x) = w_2(x)\, dx$ for which the weight
functions are related to modified Bessel functions.

The structure of the paper is as follows. In Section \ref{sec2} we will give the necessary background about
multiple orthogonal polynomials and their relation with Hermite-Pad\'e approximation.
In Section \ref{sec3} we explain what we mean by simultaneous Gaussian quadrature. 
Multiple orthogonal polynomials satisfy various recurrence relations which we give in Section \ref{sec4}, and in
particular the multiple orthogonal polynomials on the stepline (near the diagonal) satisfy a linear recurrence relation
of order $r+1$, which gives rise to a banded Hessenberg matrix with one diagonal above the main diagonal and $r$
diagonals below the main diagonal. The eigenvalues of this Hessenberg matrix correspond to the zeros of the type II
multiple orthogonal polynomial.  The main theorem (Theorem \ref{thm3})
will be proved in Section \ref{sec5}. We will use the fact that the quadrature weights can be expressed in terms 
of the Christoffel-Darboux kernel for multiple orthogonal polynomials. One important difference is that
the Christoffel-Darboux kernel for multiple orthogonal polynomials contains both the type I and the type II
multiple orthogonal polynomials, and hence is not symmetric in the two variables.
 Two numerical examples involving modified Bessel function weights are
worked out in Section \ref{sec6}.

\section{Multiple orthogonal polynomials}    \label{sec2}
Multiple orthogonal polynomials are polynomials in one variable with orthogonality relations with respect to $r$ measures
$(\mu_1,\mu_2,\ldots,\mu_r)$ (on the real line in this paper).
We use a multi-index $\vec{n}=(n_1,n_2,\ldots,n_r) \in \mathbb{N}^r$ of size $|\vec{n}|=n_1+n_2+\ldots+n_r$.

\begin{definition}[type II MOP]
A type II multiple orthogonal polynomial $P_{\vec{n}}$ is a polynomial of degree $\leq |\vec{n}|$ for which
\begin{eqnarray*}
     \int_\mathbb{R} P_{\vec{n}}(x) x^k\, d\mu_1(x) &=& 0, \qquad 0 \leq k \leq n_1-1, \\
                        & \vdots&   \\
    \int_\mathbb{R} P_{\vec{n}}(x) x^k\, d\mu_r(x) &=& 0, \qquad 0 \leq k \leq n_r-1. 
\end{eqnarray*}    
\end{definition}

\begin{definition}[type I MOP]
Type I multiple orthogonal polynomials are in a vector $(A_{\vec{n},1},\ldots,A_{\vec{n},r})$, with
$\deg A_{\vec{n},j} \leq n_j-1$, and they satisfy
\[      \sum_{j=1}^r  \int_\mathbb{R} A_{\vec{n},j}(x) x^k\, d\mu_j(x) = 0, \qquad 0 \leq k \leq |\vec {n}|-2.  \]
 \end{definition}
 
 The orthogonality conditions for the type II multiple orthogonal polynomial give a homogeneous system of $|\vec{n}|$ linear equations for the $|\vec{n}|+1$ unknown coefficients  of the polynomial $P_{\vec{n}}$. In a similar way, the orthogonality conditions for the type I multiple orthogonal polynomials give a homogeneous linear system of $|\vec{n}|-1$
equations for the $|\vec{n}|$ unknown coefficients of the polynomials $A_{\vec{n},j}$ $(1 \leq j \leq r)$.
When the type I and type II multiple orthogonal polynomials are unique up to a multiplicative factor, then we say that
the multi-index $\vec{n}$ is normal. When this happens, the type II polynomial $P_{\vec{n}}$  has degree $|\vec{n}|$ 
and we will take the leading coefficient to be one. For a normal multi-index we normalize the type I polynomials by setting
\begin{equation}   \label{typeInorm}
   \sum_{j=1}^r \int_\mathbb{R}  x^{|\vec{n}|-1} A_{\vec{n},j}(x)\, d\mu_j(x) = 1.   
\end{equation} 
In this paper we assume this is the case for all multi-indices $\vec{n}$, in which case we deal with a perfect system.
Hence our type II multiple orthogonal polynomials are monic and the type I multiple orthogonal polynomials satisfy
\eqref{typeInorm}.   
\medskip

Multiple orthogonal polynomials appear in simultaneous rational approximation of $r$ functions $(f_1,f_2,\ldots,f_r)$, where
\[          f_j(z) = \int_\mathbb{R} \frac{d\mu_j(x)}{z-x}   , \qquad   z \in \mathbb{C} \setminus \mathbb{R}.  \]

\begin{problem}[Type II Hermite-Pad\'e]
Find a polynomial $P_{\vec{n}}$ of degree $\leq |\vec{n}|$ and polynomials $Q_{\vec{n},j}$ such that for $1 \leq j \leq r$
\[     P_{\vec{n}}(z) f_j(z) - Q_{\vec{n},j}(z) = \mathcal{O}\left( \frac{1}{z^{n_j+1}} \right), \qquad z \to \infty. \]
\end{problem}

The common denominator $P_{\vec{n}}$ is the type II multiple orthogonal polynomial and the numerator polynomials
are given by
\[    Q_{\vec{n},j}(z) = \int_\mathbb{R} \frac{P_{\vec{n}}(z)-P_{\vec{n}}(x)}{z-x}\, d\mu_j(x).  \]
The error can also be expressed in terms of the type II multiple orthogonal polynomial:
\[    P_{\vec{n}}(z) f_j(z) - Q_{\vec{n},j}(z) = \int_\mathbb{R} \frac{P_{\vec{n}}(x)}{z-x}\, d\mu_j(x).   \]

\begin{problem}[Type I Hermite-Pad\'e]
Find polynomials $A_{\vec{n},j}$ of degree $\leq n_j-1$ and a polynomial $B_{\vec{n}}$ such that 
\[    \sum_{j=1}^r A_{\vec{n},j}(z) f_j(z) - B_{\vec{n}}(z) = \mathcal{O}\left( \frac{1}{z^{|\vec{n}|}} \right), \qquad z \to \infty. \]
\end{problem}

The polynomials $A_{\vec{n},1},\ldots,A_{\vec{n},r}$ are the type I multiple orthogonal polynomials and the polynomial
$B_{\vec{n}}$ is
\[     B_{\vec{n}}(z) = \sum_{j=1}^r \int_\mathbb{R} \frac{A_{\vec{n},j}(z)-A_{\vec{n},j}(x)}{z-x}\, d\mu_j(x).  \]
The error is given by
 \[    \sum_{j=1}^r A_{\vec{n},j}(z) f_j(z) - B_{\vec{n}}(z) 
 = \sum_{j=1}^r \int_\mathbb{R} \frac{A_{\vec{n},j}(x)}{z-x}\, d\mu_j(x).  \]

Here is some history indicating important progress (but very incomplete).
Hermite-Pad\'e approximation was introduced by Charles Hermite in 1873, who used them in his proof that $e$ is transcendental \cite{Hermite}. The case $r=1$ was then investigated in more detail by Hermite's student Henri Pad\'e (Pad\'e approximation, Pad\'e table) \cite{Pade}. One important system for which all the multi-indices are normal (a perfect system)
was described by Aurel Angelescu \cite{Angelesco} who used $r$ measures (weights) supported on $r$ disjoint intervals.
This is nowadays known as an Angelesco system. In 1934 Kurt Mahler \cite{Mahler} investigated perfect systems but his work was published much later 1968. Another important system of measures was introduced by Evgeni Nikishin in 1979 \cite{Nikishin}. 
For $r=2$ the measures $(\mu_1,\mu_2)$ are absolutely continuous on an interval $\Delta_0$ with weight functions
$(w_1,w_2)$ for which 
\[   \frac{w_2(x)}{w_1(x)} = \int_{\Delta_1} \frac{d\sigma(t)}{x-t}, \]
where $\sigma$ is a positive measure on an interval $\Delta_1$ which is disjoint from $\Delta_0$. 
Such a system (and its recursive extensions to $r >2$) are now known as Nikishin systems and it has been shown that
these are perfect systems \cite{FidalgoLopez}.
Gonchar and Rakhmanov studied the convergence of Hermite-Pad\'e approximants in \cite{GoncharRakh} for which they
used an extremal problem in logarithmic potential theory for a vector of measures.
The idea of  simultaneous Gaussian quadrature was introduced by Carlos F. Borges in \cite{Borges}. 
Van Assche, Geronimo and Kuij\-laars formulated a Riemann-Hilbert problem for multiple orthogonal polynomials in 2001 \cite{WVAGerKuijl}. This allows to obtain the asymptotic behavior of multiple orthogonal polynomials.

Nowadays multiple orthogonal polynomials are used in various applications, such as eigenvalues of
products of random matrices, non-intersecting Brownian motions, combinatorial problems related to
random tilings, rational approximations in number theory, etc. They are also used in numerical analysis for matrix function evaluation (Alqahtani and Reichel \cite{Alqahtani}).

\subsection{Examples extending classical orthogonal polynomials}

Many examples of multiple orthogonal polynomials are known. 
\begin{itemize}
\item Multiple Hermite polynomials: the weights are normal densities with non-zero mean
\[    w_j(x) = e^{-x^2+c_jx}, \quad c_i \neq c_j,  \qquad x \in (-\infty,\infty). \]
\item Multiple Laguerre polynomials: the weights
\[     w_j(x) = x^{\alpha_j} e^{-x}, \quad  \alpha_i-\alpha_j \notin \mathbb{Z},\  \alpha_j >-1,  \qquad x \in [0,\infty) \]
correspond to multiple Laguerre polynomials of the first kind. The weights
\[     w_j(x) = x^\alpha e^{-c_jx}, \quad c_i \neq c_j, \ c_j > 0,\ \alpha >-1, \qquad x \in [0,\infty), \]
correspond to multiple Laguerre polynomials of the second kind.
\item Jacobi-Pi\~neiro polynomials: the weights are
\[    w_j(x) = x^{\alpha_j} (1-x)^\beta, \quad    \alpha_i-\alpha_j \notin \mathbb{Z},\  \alpha_j,\beta > -1, \qquad x \in [0,1].  \]
\end{itemize}
For these extensions of the classical orthogonal polynomials one knows a differential equation for the multiple orthogonal polynomials, a system of nearest neighbor recurrence relations, the asymptotic distribution of the zeros, etc.
Other examples include weights in terms of special functions satisfying a second order differential
equation.
\subsection{Modified Bessel functions $(K_\nu,K_{\nu+1})$}
The modified Bessel function $K_\nu$ is given as
\[   K_\nu(x) = \frac12 \left( \frac{x}{2}\right)^\nu \int_0^\infty \exp \left(-t-\frac{x^2}{4t} \right) t^{-\nu-1}\, dt. \]
It  satisfies the modified Bessel equation
\[  x^2 y'' + x y' - (x^2+ \nu^2) y = 0.  \]
Introduce the functions
\[    \rho_\nu(x) = 2 x^{\nu/2} K_\nu(2\sqrt{x}), \qquad x >0, \]
which have the moments
\[   \int_0^\infty x^n \rho_\nu(x)\, dx = \Gamma(n+\nu+1) \Gamma(n+1)   . \]
Van Assche and Yakubovich \cite{WVAYakubovich} and Ben Cheikh and Douak \cite{BenCheikhDouak}
obtained the multiple orthogonal polynomials for the weights or $(w_1,w_2)=x^\alpha (\rho_\nu, \rho_{\nu+1})$,
where $\alpha > -1$ and $\nu \geq 0$. The vector Pearson equation is 
\[   \left[ x \begin{pmatrix} w_1 \\ w_2 \end{pmatrix}\right]' 
 =\begin{pmatrix} \alpha + \nu+1 & -1 \\ -x & \alpha+1 \end{pmatrix}
    \begin{pmatrix} w_1 \\ w_2 \end{pmatrix}  . \]
One has
\[   \frac{\rho_\nu(x)}{\rho_{\nu+1}(x)} = \frac{1}{\pi^2} \int_0^\infty 
     \frac{1}{(x+s) [J_{\nu+1}^2(2\sqrt{s})+ Y_{\nu+1}^2(2\sqrt{s})]} \frac{ds}{s}, \]
so that this is a Nikishin system.
Denote the type I function by
\[  q_{n,m}^\alpha(x) = A_{n,m}(x) \rho_\nu(x) + B_{n,m}(x) \rho_{\nu+1}(x),   \]
and the type II multiple orthogonal polynomials by $p_{n,m}^{\alpha}$.
\begin{property}[Rodrigues formula]    \label{prop1}
One has
\[  x^\alpha q_{n,n}^\alpha = \frac{d^{2n-1}}{dx^{2n-1}} [x^{2n-1+\alpha} \rho_\nu(x)], \]
\[    x^\alpha q_{n,n-1}^\alpha = \frac{d^{2n-2}}{dx^{2n-2}} [x^{2n+\alpha-2} \rho_\nu(x)].   \]
\end{property}    
\begin{property}[recurrence relation]    \label{prop2}
Let $P_{2n}=p_{n,n}^\alpha$ and $P_{2n+1}=p_{n+1,n}^{\alpha}$, then
\[  xP_n(x) = P_{n+1}(x) + b_nP_n(x) + c_nP_{n-1}(x) +d_nP_{n-2}(x)  \]
with
\begin{eqnarray*}
   b_n &=& (n+\alpha+1)(3n+\alpha+2\nu)-(\alpha+1)(\nu-1), \\
   c_n &=& n(n+\alpha)(n+\alpha+\nu)(3n+2\alpha+\nu) ,  \\
   d_n &=& n(n-1)(n+\alpha)(n+\alpha-1)(n+\alpha+\nu)(n+\alpha+\nu-1).  
\end{eqnarray*}   
Let $Q_{2n}=q_{n,n}^\alpha$ and $Q_{2n-1}=q_{n,n-1}^\alpha$, then
\[   xQ_n(x) = Q_{n-1}(x) + b_{n-1} Q_n(x) + c_n Q_{n+1}(x) + d_{n+1} Q_{n+2}(x). \]
\end{property}

\subsection{Modified Bessel functions $(I_{\nu},I_{\nu+1})$}
Another solution of the modified Bessel equation is given by the modified Bessel function $I_\nu(x)$.
In \cite{ECoussWVA} we considered
\[  w_1(x) = x^{\nu/2}I_\nu(2\sqrt{x}) e^{-cx}, \quad   w_2(x) = x^{(\nu+1)/2} I_{\nu+1}(2\sqrt{x}) e^{-cx}, 
\quad x>0, \]
with $\nu>-1$ and $c>0$. The case $c=1$ was first studied in \cite{Douak}.
These weights are related to the noncentral $\chi^2$-distribution:
\[                   p_{\chi^2}^{(\nu,\lambda)}(x) = \frac12 \left( \frac{x}{\lambda} \right)^{(\nu-2)/4} I_{(\nu-2)/2}(\sqrt{\lambda x}) e^{-(\lambda+x)/2}, \qquad x > 0,  \]
with $\lambda > 0$ and $\nu \in \mathbb{N}$.
The Pearson equation is
\[  \left[ x \begin{pmatrix} w_1 \\ w_2 \end{pmatrix} \right]' = \begin{pmatrix}\nu-cx & 1 \\ x & -cx \end{pmatrix}
   \begin{pmatrix} w_1 \\ w_2  \end{pmatrix}  .  \]
One has
\[   \frac{w_1(x)}{w_2(x)} = \frac{\nu+1}{x} +  \sum_{n=1}^\infty \frac{1}{x+j_{\nu+1,n}^2/4}
     \frac{J_{\nu+2}(j_{\nu+1,n})}{J_{\nu+1}'(j_{\nu+1,n})}, \]
so that we again have a Nikishin system. 
Let
\[ Q_{2n}^{\nu}(x) = q_{n,n}^{\nu}(x), \quad 
   Q_{2n+1}^{\nu}(x)=q_{n+1,n}^{\nu}(x), \]
       and similarly 
\[     P_{2n}^{\nu}(x) = p_{n,n}^\nu(x), \quad P_{2n+1}^\nu(x)=p_{n+1,n}^\nu(x).  \]
\begin{property}[raising-lowering]   \label{prop3}
For every $\nu>-1$ and $c>0$ one has
\[   \frac{d}{dx} Q_n^{\nu+1}(x) = Q_{n+1}^\nu(x),  \quad \textup{ and } \quad
    \frac{d}{dx} P_n^\nu(x) = n P_{n-1}^{\nu+1}(x).  \]
\end{property}
\begin{property}[recurrence relation]    \label{prop4}
One has
\[   x P_n^\nu(x) = P_{n+1}^\nu(x) + b_n P_n^\nu(x) + c_n P_{n-1}^\nu(x) + d_n P_{n-2}^\nu(x), \]
with
\begin{eqnarray*}
  b_n &=& \frac{1}{c^2} [1+c(\nu+2n+1)], \\
  c_n &=& \frac{n}{c^3} [2+c(\nu+n)], \\
  d_n &=& \frac{n(n-1)}{c^4}.
\end{eqnarray*}
Furthermore
\[  xQ_n^\nu(x) = Q_{n-1}^\nu(x) + b_{n-1} Q_n^\nu(x) + c_n Q_{n+1}^\nu(x) + d_{n+1} Q_{n+2}^\nu(x).   \]
\end{property}

\section{Simultaneous Gaussian quadrature}  \label{sec3}
Suppose we have two measures $\mu_1,\mu_2$ and one function $f$. We want to approximate integrals by sums
\[     \int_\mathbb{R} f(x) \, d\mu_1(x)   \approx   \sum_{k=1}^{2n} \lambda_{k,2n}^{(1)} f(x_{k,2n}), \]
\[     \int_\mathbb{R} f(x) \, d\mu_2(x)   \approx   \sum_{k=1}^{2n} \lambda_{k,2n}^{(2)} f(x_{k,2n}). \]

\begin{theorem}[Borges \cite{Borges}]
If we take for $\{x_{k,2n}, 1 \leq k \leq 2n \}$ the zeros of the type II multiple orthogonal $P_{n,n}$
for the two measures $(\mu_1,\mu_2)$ and use interpolatory quadrature, then the quadrature  is
exact for polynomials $f$ of degree $\leq 3n-1$.
\end{theorem}

There are $2n$ function evaluations. If we use ordinary Gaussian quadrature with $n$ nodes for
each integral, then we also use $2n$ function evaluations, but the quadrature then is exact for
polynomials $f$ of degree $\leq 2n-1$. Important numerical aspects are:
\begin{itemize}
\item How to compute the quadrature nodes $x_{k,2n}$, $1 \leq k \leq 2n$ (the zeros of $P_{n,n}$)?
\item How to compute the quadrature weights
\begin{equation} \label{lambda1}
     \lambda_{k,2n}^{(1)} = \int_\mathbb{R} \frac{P_{n,n}(x)}{(x-x_{k,2n}) P_{n,n}'(x_{k,2n})} \, d\mu_1(x) , 
\end{equation}     
\begin{equation} \label{lambda2}
     \lambda_{k,2n}^{(2)} = \int_\mathbb{R} \frac{P_{n,n}(x)}{(x-x_{k,2n}) P_{n,n}'(x_{k,2n})} \, d\mu_2(x) .
\end{equation}     
\end{itemize}
Some relevant analytical problems are
\begin{itemize}
\item The convergence of the quadrature rules as $n \to \infty$.
\item Are the weights $\lambda_{k,2n}^{(1)}$ and $\lambda_{k,2n}^{(2)}$ positive? If not, are they small?
\end{itemize}
We have investigated some of the analytical problems of simultaneous Gaussian quadrature for an Angelesco system in \cite{LubWVA} and for multiple Hermite polynomials in \cite{WVAVuer}.

We formulated the simultaneous Gaussian quadrature for $2n$ nodes. If the number of nodes is odd ($2n+1$), then
one uses the zeros of the type II multiple orthogonal $P_{n+1,n}$ and
\[    \int_\mathbb{R} f(x) \, d\mu_1(x)   =   \sum_{k=1}^{2n+1} \lambda_{k,2n+1}^{(1)} f(x_{k,2n+1})  \]
holds for polynomials $f$ of degree $\leq 3n+1$, and
\[    \int_\mathbb{R} f(x) \, d\mu_2(x)   =   \sum_{k=1}^{2n+1} \lambda_{k,2n+1}^{(2)} f(x_{k,2n+1})  \]
holds for polynomials $f$ of degree $\leq 3n$. For multi-indices $(n,m)$ one takes the nodes at the zeros of $P_{n,m}$
and
\[    \int_\mathbb{R} f(x) \, d\mu_1(x)   =   \sum_{k=1}^{n+m} \lambda_{k,n+m}^{(1)} f(x_{k,n+m}), \qquad \deg f \leq 2n+m-1,  \]
\[    \int_\mathbb{R} f(x) \, d\mu_2(x)   =   \sum_{k=1}^{n+m} \lambda_{k,n+m}^{(2)} f(x_{k,n+m}), \qquad \deg f \leq n+2m-1.  \]

There is an important relation with Hermite-Pad\'e approximation.
Recall that for type II Hermite-Pad\'e approximation one has
\[     P_{\vec{n}}(z) f_j(z) - Q_{\vec{n},j}(z) = \mathcal{O}\left( \frac{1}{z^{n_j+1}} \right), \qquad z \to \infty. \]
For the rational approximants one then has
\[     f_j(z) - \frac{Q_{\vec{n},j}(z)}{P_{\vec{n}}(z)} = \mathcal{O} \left( \frac{1}{z^{|\vec{n}|+n_j+1}} \right), \]
and if one decomposes the rational approximant into partial fractions, then 
\[  \frac{Q_{\vec{n},j}(z)}{P_{\vec{n}}(z)} = 
\sum_{k=1}^{|\vec{n}|} \frac{\lambda_{k,\vec{n}}^{(j)}}{z-x_{k,\vec{n}}} .   \]
Hence the poles of the type II Hermite-Pad\'e approximants are the quadrature nodes and the residues are the quadrature
weights. This extends the well-known connection between Gaussian quadrature and Pad\'e approximation.
If $\mu_j$ is supported on $[a,b]$ and $g$ is analytic in a neighborhood $\Omega$ of $[a,b]$, then the error is given by
\[   \int_a^b g(x) \, d\mu_j(x) - \sum_{k=1}^{|\vec{n}|} \lambda_{k,\vec{n}}^{(j)} g(x_{k,\vec{n}}) 
   = \frac{1}{2\pi i} \int_\Gamma g(z) \left( f_j(z) 
      -  \frac{Q_{\vec{n},j}(z)}{P_{\vec{n}}(z)} \right) \, dz  ,   \]
 where $\Gamma$ is a closed contour in $\Omega$ around $[a,b]$. Hence results about the asymptotic behavior of 
 Hermite-Pad\'e approximation imply results for the convergence of the quadrature formula when the integrand 
$g$ is sufficiently smooth.
 
\section{Recurrence relations}     \label{sec4}
Denote by $\vec{e}_j = (0,\ldots,0,1,0,\ldots,0)$ the standard unit vectors. 
Multiple orthogonal polynomials always satisfy a system of linear recurrence relations connecting $P_{\vec{n}}$
to its nearest neighbors $P_{\vec{n}+\vec{e}_k}$ from above and $P_{\vec{n}-\vec{e}_j}$ from below
\cite{WVA}:
\begin{equation}   \label{nearrecP}
  xP_{\vec{n}}(x) = P_{\vec{n}+\vec{e}_k}(x) + b_{\vec{n},k} P_{\vec{n}}(x)
   + \sum_{j=1}^r a_{\vec{n},j} P_{\vec{n}-\vec{e}_j}(x), \qquad 1 \leq k \leq r. 
\end{equation}
For the type I functions 
\[  Q_{\vec{n}}(x) = \sum_{j=1}^r  A_{\vec{n},j}(x) w_j(x), \qquad   w_j(x) = \frac{d\mu_j(x)}{d\mu(x)}, \]
where $\mu=\mu_1+\mu_2+\cdots+\mu_r$, the nearest neighbor recurrence relations are
\begin{equation}   \label{nearrecQ}
  xQ_{\vec{n}}(x) = Q_{\vec{n}-\vec{e}_k}(x) + b_{\vec{n}-\vec{e}_k,k} Q_{\vec{n}}(x)
   + \sum_{j=1}^r a_{\vec{n},j} Q_{\vec{n}+\vec{e}_j}(x), \qquad 1 \leq k \leq r. 
\end{equation}
If we combine these $r$ recurrence relations, then one can find a $(r+2)$-term recurrence relation on
 the stepline: for $n=kr+\ell$ we take $\vec{n}=(\underbrace{k+1,\ldots,k+1}_{\text{$\ell$ times}},k,\ldots,k)$ and 
 we define the type II polynomials on the stepline by $P_n(x) = P_{\vec{n}}(x)$. This $(r+2)$-term recurrence relation is 
\begin{equation}  \label{steprec}
   xP_n(x) = P_{n+1}(x) + \sum_{j=0}^r c_{n,j} P_{n-j}(x).   
\end{equation}   
There is also a $(r+2)$-term recurrence relation for the type I functions on the stepline:
\begin{equation}  \label{steprecQ}
   xQ_n(x) = Q_{n-1}(x) + \sum_{j=0}^r c_{n+j-1,j} Q_{n+j}(x).   
\end{equation} 

Filipuk, Haneczok and Van Assche \cite{FilHanWVA} worked out an algorithm to go from the nearest neighbor recurrence coefficients  $(a_{\vec{n},j},b_{\vec{n},j};1 \leq j \leq r)_{\vec{n}\in \mathbb{N}^r}$ to the stepline
coefficients  $(c_{n,0},c_{n,1},\ldots,c_{n,r})_{n \in \mathbb{N}}$. They also gave algorithms to go from
the recurrence coefficients of the orthogonal polynomials of the measures $\mu_j$ to the nearest neighbor recurrence
coefficients. We will illustrate this for the case $r=2$.
The nearest neighbor recurrence relations are
\begin{eqnarray}
    xP_{n,m}(x) &=& P_{n+1,m}(x) + c_{n,m}P_{n,m}(x) + a_{n,m}P_{n-1,m}(x) + b_{n,m} P_{n,m-1}(x) , \label{NNP1} \\
    xP_{n,m}(x) &=& P_{n,m+1}(x) + d_{n,m}P_{n,m}(x) + a_{n,m}P_{n-1,m}(x) + b_{n,m} P_{n,m-1}(x) . \label{NNP2}
\end{eqnarray}
Subtracting both equations gives
\begin{equation}  \label{Psub}
    P_{n+1,m}(x) - P_{n,m+1}(x) = (d_{n,m}-c_{n,m}) P_{n,m}(x).
\end{equation} 
Denote the type II polynomials on the stepline by $P_{2k}(x)=P_{k,k}(x)$, $P_{2k+1}(x)=P_{k+1,k}(x)$.
Take $m=n$ then \eqref{NNP1} is
\[  x P_{n,n}(x) = P_{n+1,n}(x) + c_{n,n} P_{n,n}(x) + a_{n,n} P_{n-1,n}(x) + b_{n,n} P_{n,n-1}(x).  \]
Use \eqref{Psub} to replace $P_{n-1,n}(x)$ by $P_{n,n-1}(x) + (c_{n-1,n-1}-d_{n-1,n-1})P_{n-1,n-1}(x)$ to find
\begin{multline*}
  xP_{2n}(x) = P_{2n+1}(x) + c_{n,n} P_{2n}(x) + (a_{n,n}+b_{n,n}) P_{2n-1}(x) \\
  + a_{n,n} (c_{n-1,n-1}-d_{n-1,n-1})    P_{2n-2}(x).  
\end{multline*}  
In a similar way, by taking $m=n-1$, \eqref{NNP2} is
\[  xP_{n,n-1}(x) = P_{n,n}(x) + d_{n,n-1}P_{n,n-1}(x) + a_{n,n-1}P_{n-1,n-1}(x) + b_{n,n-1} P_{n,n-2}(x).  \]
Use \eqref{Psub} to replace $P_{n,n-2}(x)$ by $P_{n-1,n-1}(x) + (d_{n-1,n-2}-c_{n-1,n-2}) P_{n-1,n-2}(x)$ to find       
\begin{multline*}
  xP_{2n-1}(x) = P_{2n}(x) + d_{n,n-1}P_{2n-1}(x) + (a_{n,n-1}+b_{n,n-1}) P_{2n-2}(x)   \\
   + b_{n,n-1}(d_{n-1,n-2}-c_{n-1,n-2}) P_{2n-3}(x).  
\end{multline*}   
This means that the four-term recurrence relation for the type II polynomials on the stepline is
\begin{equation}  \label{Prec}
   xP_n(x) = P_{n+1}(x) + b_n P_n(x) + c_n P_{n-1}(x) + d_n P_{n-2}(x),
\end{equation}   
with
\begin{align*}
    b_{2k} &= c_{k,k}, \quad b_{2k+1} = d_{k+1,k} \\
    c_{2k} &= a_{k,k}+b_{k,k}, \quad c_{2k+1} = a_{k+1,k}+b_{k+1,k}, \\
    d_{2k} &= a_{k,k}(c_{k-1,k-1}-d_{k-1,k-1}), \quad d_{2k+1} = b_{k+1,k}(d_{k,k-1}-c_{k,k-1}).
     \end{align*}
For the type I functions, the nearest neighbor recurrence relations are
\begin{eqnarray*}
  xQ_{n,m}(x) &=& Q_{n-1,m}(x) + c_{n-1,m} Q_{n,m}(x) + a_{n,m} Q_{n+1,m}(x) + b_{n,m}(x) Q_{n,m+1}(x), \\
  xQ_{n,m}(x) &=& Q_{n,m-1}(x) + d_{n,m-1} Q_{n,m}(x) + a_{n,m} Q_{n+1,m}(x) + b_{n,m}(x) Q_{n,m+1}(x).
\end{eqnarray*} 
Subtracting both equations gives
\begin{equation}  \label{Qsub}
   Q_{n-1,m}(x)-Q_{n,m-1}(x) = (d_{n,m-1}-c_{n-1,m}) Q_{n,m}(x).
\end{equation}
In a similar way as before, one then finds for the type I functions on the stepline $Q_{2k}(x)=Q_{k,k}(x)$, $Q_{2k+1}(x)=Q_{k+1,k}(x)$
 \begin{equation}  \label{Qrec}
    xQ_n(x) = Q_{n-1}(x) + b_{n-1} Q_n(x) + c_n Q_{n+1}(x) + d_{n+1} Q_{n+2}(x) .
\end{equation}        

Using the coefficients in this four-term recurrence relation, we can create a banded Hessenberg matrix
\begin{equation}      \label{Hessenberg}
    H_N = \begin{pmatrix}
         b_0 & 1 & 0 & 0 & 0 & 0  & \cdots & 0 \\
         c_1 & b_1 & 1 & 0 & 0 & 0 & \cdots & 0 \\
         d_2 & c_2 & b_2 & 1 & 0 & 0 & \cdots & 0 \\
          0 & d_3 & c_3 & b_3 & 1 & 0 & \cdots & 0 \\
          0 & 0 & d_4 & c_4 & b_4 & 1 & \cdots & 0 \\
          \vdots &  &   & \ddots & \ddots & \ddots & \ddots & \vdots \\
          0 & \cdots & 0 & 0 & d_{N-2} & c_{N-2} & b_{N-2} & 1 \\
          0 & \cdots& 0 & 0 & 0& d_{N-1} & c_{N-1} & b_{N-1}
          \end{pmatrix}  .
\end{equation}          
From now one we will use $N$ for the number of quadrature nodes. If $N=2n$ then this corresponds to the multi-index $(n,n)$
and for $N=2n+1$ to the multi-index $(n+1,n)$.
The recurrence relation \eqref{Prec} for the type II polynomials can then be written as
\[   H_N \begin{pmatrix} P_0(x) \\ P_1(x) \\ \vdots \\ P_{N-1}(x) \end{pmatrix}
    + \begin{pmatrix} 0 \\  \vdots \\ 0  \\ P_N(x) \end{pmatrix} 
 =  x \begin{pmatrix} P_0(x) \\ P_1(x) \\ \vdots \\ P_{N-1}(x) \end{pmatrix} .
   \]
If we evaluate this at a zero $x_{k,N}$ of $P_N$, then we get the following result:   
\begin{property}  \label{prop5}
If $P_N(x_{k,N}) = 0$ then $x_{k,N}$ is an eigenvalue of $H_N$ with \textbf{right} eigenvector 
$\bigl(P_0(x_{k,N}),P_1(x_{k,N}),\ldots,P_{N-1}(x_{k,N})\bigr)$.
\end{property}

There is a similar property for the type I multiple orthogonal polynomials $(A_N,B_N)$ for which
\[    A_{2k}(x) = A_{k,k}(x), \quad A_{2k+1}(x) = A_{k+1,k}(x), \]
\[    B_{2k}(x) = B_{k,k}(x), \quad B_{2k+1}(x) = B_{k+1,k}(x). \]

\begin{property}   \label{prop6}
If $P_N(x_{k,N}) = 0$ then $x_{k,N}$ is an eigenvalue of $H_N$ with \textbf{left} eigenvector 
$\bigl(C_1A_1(x_{k,N})+C_2B_1(x_{k,N}),C_1A_2(x_{k,N})+C_2B_2(x_{k,N}),\ldots,C_1A_n(x_{k,N})+C_2B_n(x_{k,N})\bigr)$,
with constants $C_1$ and $C_2$ satisfying 
\begin{equation}   \label{C1C2}
  \begin{cases}
      C_1A_{N+1}(x_{k,N}) + C_2B_{N+1}(x_{k,N}) = 0, \\
      C_1A_{N+2}(x_{k,N}) + C_2 B_{N+2}(x_{k,N}) = 0.
      \end{cases} 
\end{equation}  
\end{property}

\begin{proof} For the type I function $Q_n(x) = A_n(x) w_1(x) + B_n(x) w_2(x)$, the four-term recurrence
relation \eqref{Qrec} is equivalent to
\begin{align*}
  &\begin{pmatrix} Q_1(x) & Q_2(x) & \cdots & Q_N(x) \end{pmatrix} H_N \\
    &\quad + \begin{pmatrix} 0 & \cdots & 0 &  d_N Q_{N+1}(x) &\qquad  0 \qquad \ \ \ \end{pmatrix} \\
     &\quad +  \begin{pmatrix} 0 & \cdots & 0 & c_N Q_{N+1}(x) & d_{N+1} Q_{N+2}(x) \end{pmatrix} \\
     &= x \begin{pmatrix} Q_1(x) & Q_2(x) & \cdots & Q_N(x) \end{pmatrix} .
\end{align*}
This relation also holds for the polynomials $A_n(x)$ and $B_n(x)$, 
hence 
\begin{align*}
  &\begin{pmatrix} A_1(x) & A_2(x) & \cdots & A_N(x) \end{pmatrix} H_N \\
    &\quad + \begin{pmatrix} 0 & \cdots & 0 &  d_N A_{N+1}(x) &\qquad  0 \qquad \ \ \ \end{pmatrix} \\
     &\quad +  \begin{pmatrix} 0 & \cdots & 0 & c_N A_{N+1}(x) & d_{N+1} A_{N+2}(x) \end{pmatrix} \\
     &= x \begin{pmatrix} A_1(x) & A_2(x) & \cdots & A_N(x) \end{pmatrix} ,
\end{align*}
and
\begin{align*}
  &\begin{pmatrix} B_1(x) & B_2(x) & \cdots & B_N(x) \end{pmatrix} H_N \\
    &\quad + \begin{pmatrix} 0 & \cdots & 0 &  d_N B_{N+1}(x) &\qquad  0 \qquad \ \ \ \end{pmatrix} \\
     &\quad +  \begin{pmatrix} 0 & \cdots & 0 & c_N B_{N+1}(x) & d_{N+1} B_{N+2}(x) \end{pmatrix} \\
     &= x \begin{pmatrix} B_1(x) & B_2(x) & \cdots & B_N(x) \end{pmatrix} ,
\end{align*}
and therefore also for a linear combination of $A_n$ and $B_n$
\begin{align*}  
   &\begin{pmatrix}C_1 A_1(x)+C_2B_1(x) & C_1A_2(x)+C_2B_2(x) & \cdots & C_1A_N(x)+C_2B_N(x) \end{pmatrix} H_N \\
    & + \begin{pmatrix} 0 & \cdots & 0 &  d_N [C_1A_{N+1}(x)+C_2B_{N+1}(x)] &\qquad \qquad \qquad \ 0 \qquad \qquad \qquad \end{pmatrix} \\
     & +  \begin{pmatrix} 0 & \cdots & 0 & c_N[C_1 A_{N+1}(x)+C_2B_{N+1}(x)] & d_{N+1} [C_1 A_{N+2}(x)+C_2B_{N+2}(x)] \end{pmatrix} \\
     &= x \begin{pmatrix} C_1A_1(x)+C_2B_1(x) & C_1A_2(x)+C_2B_2(x) & \cdots & C_1A_N(x)+C_2B_N(x) \end{pmatrix} .
\end{align*}
Choose $C_1$ and $C_2$ such that
\begin{equation*}   
  \begin{cases}
      C_1A_{N+1}(x) + C_2B_{N+1}(x) = 0, \\
      C_1A_{N+2}(x) + C_2 B_{N+2}(x) = 0.
      \end{cases} 
\end{equation*}   
This is a homogeneous system of equations for $C_1$ and $C_2$, which only has a non-trivial
solution if the matrix is singular. Mahler's relation between type I and type II polynomials \cite[Eq. (4.2)]{WVAGerKuijl}
gives
\[      \det \begin{pmatrix}  A_{N+1}(x) & B_{N+1}(x) \\
                                          A_{N+2}(x) & B_{N+2}(x)  \end{pmatrix} = C_3 P_N(x)  \]
for some constant $C_3$. Hence if $x_{k,N}$ is a zero of $P_N$, then we can find $(C_1,C_2)\neq(0,0)$  so that \eqref{C1C2} holds.
\end{proof}
Observe that $C_1$ and $C_2$ may depend on the eigenvalue $x_{k,N}$.

\section{Golub-Welsch type result}    \label{sec5}
We are now ready to prove the main result. 
\begin{theorem}    \label {thm3}
The quadrature nodes $\{x_{j,N}, 1 \leq j \leq N\}$ for simultaneous Gaussian quadrature
are given by the eigenvalues of the banded Hessenberg matrix $H_N$.
Let $\vec{v}_{N,j}$ be the right eigenvector of $H_N$ for the eigenvalue $x_{j,N}$ with first component $1$ 
and $\vec{u}_{N,j}$ be a left eigenvector for $x_{j,N}$. Then
\begin{eqnarray*}   
\lambda_{j,N}^{(1)} &=& \frac{D_{1,1} \vec{u}_{N,j}(1)}{\langle \vec{u}_{N,j}, \vec{v}_{N,j} \rangle},  \\
\lambda_{j,N}^{(2)} &=& \frac{D_{2,1} \vec{u}_{N,j}(1) + D_{2,2} \vec{u}_{N,j}(2)}{\langle \vec{u}_{N,j}, \vec{v}_{N,j} \rangle},  
\end{eqnarray*}
with constants $D_{i,j}$ given by
\[   \begin{pmatrix}  D_{1,1} & 0 \\ D_{2,1} & D_{2,2} \end{pmatrix}
     = \begin{pmatrix}  A_1 & 0 \\  A_2 & B_2   \end{pmatrix}^{-1}.   \]
\end{theorem}     
This result was given for general $r$ in Coussement-Van Assche \cite[Thm.3.2]{JCoussWVA}.

\begin{proof}
We first need the Christoffel-Darboux formula for multiple orthogonal polynomials (Daems-Kuijlaars \cite{DaemsKuijl}).
Let $N=|\vec{n}|$ and $(\vec{n}_k)_{0 \leq k \leq N}$ be a path in $\mathbb{N}^r$ from
$\vec{n}_0=\vec{0}$ to $\vec{n}_N = \vec{n}$ such that $\vec{n}_{k+1}-\vec{n}_k=\vec{e}_j$ for some $1 \leq j \leq r$. Then
\[  (x-y) \sum_{k=0}^{N-1} P_{\vec{n}_k}(x) Q_{\vec{n}_{k+1}}(y)
    = P_{\vec{n}}(x)Q_{\vec{n}}(y) - \sum_{j=1}^r a_{\vec{n},j} P_{\vec{n}-\vec{e}_j}(x)
        Q_{\vec{n}+\vec{e}_j}(y).  \]
Here we used the version as it is given in terms of the nearest neighbor recurrence coefficients \cite{WVA}.
For $r=2$ and on the stepline, it becomes
\begin{multline*}   
(x-y) \sum_{k=0}^{N-1} P_k(x)Q_{k+1}(y) = P_{N}(x)Q_{N}(y) - c_{N}P_{N-1}(x)Q_{N+1}(y)  \\
- d_{N}P_{N-2}(x)Q_{N+1}(y) - d_{N+1}P_{N-1}(x)Q_{N+2}(y).
\end{multline*}
The Christoffel-Darboux formula also holds if we replace $Q_k(y)$ by $C_1A_k(y) + C_2B_k(y)$:
\begin{eqnarray*}   
\lefteqn{(x-y) \sum_{k=0}^{N-1} P_k(x)[C_1A_{k+1}(y)+C_2B_{k+1}(y)]} & & \\
  &=&   P_{N}(x)[C_1A_{N}(y)+C_2B_{N}(y)] - c_{N}P_{N-1}(x)[C_1A_{N+1}(y)+C_2B_{N+1}(y)]  \\
   & & -\ d_{N}P_{N-2}(x)[C_1A_{N+1}(y)+C_2B_{N+1}(y)] - d_{N+1}P_{N-1}(x)[C_1A_{N+2}(y)+C_2B_{N+2}(y)].
\end{eqnarray*}
Take $y=x_{j,N}$ a zero of $P_{N}$, then \eqref{C1C2} gives
\begin{equation}   \label{CDx-xj}
 \sum_{k=0}^{N-1} P_k(x)[C_1A_{k+1}(x_{j,N})+C_2B_{k+1}(x_{j,N})]  \\
=    \frac{P_{N}(x)}{x-x_{j,N}} [C_1A_{N}(x_{j,N})+C_2B_{N}(x_{j,N})].    
\end{equation} 
Recall that the quadrature nodes are given by \eqref{lambda1} and \eqref{lambda2}.
Integrating \eqref{CDx-xj} with measure $\mu_1$ gives
\begin{multline*}
      \sum_{k=0}^{N-1} [C_1A_{k+1}(x_{j,N})+C_2B_{k+1}(x_{j,N})] \int_\mathbb{R} P_k(x) \, d\mu_1(x) \\
     = \int_\mathbb{R} \frac{P_{N}(x)}{x-x_{j,N}} \, d\mu_1(x) \  [C_1A_{N}(x_{j,N})+C_2B_{N}(x_{j,N})] ,
\end{multline*}     
hence if we retain only the non-vanishing integrals
\begin{equation}  \label{lambda1a}
    \lambda_{j,N}^{(1)}  [C_1A_{N}(x_{j,N})+C_2B_{N}(x_{j,N})]  =
     \frac{1}{P_{N}'(x_{j,N})} \int_\mathbb{R} P_0(x) \, d\mu_1(x) \  [C_1A_1 + C_2 B_1]   .  
\end{equation}
In a similar way, integrating \eqref{CDx-xj} with measure $\mu_2$ gives
\begin{multline*}
      \sum_{k=0}^{N-1} [C_1A_{k+1}(x_{j,N})+C_2B_{k+1}(x_{j,N})] \int_\mathbb{R} P_k(x) \, d\mu_2(x) \\
     = \int_\mathbb{R} \frac{P_{N}(x)}{x-x_{j,N}} \, d\mu_2(x) \ [C_1A_{N}(x_{j,N})+C_2B_{N}(x_{j,N})],  
\end{multline*}  
hence if we retain only the non-vanishing integrals
\begin{multline}   \label{lambda2a}
    \lambda_{j,N}^{(2)}  [C_1A_{N}(x_{j,N})+C_2B_{N}(x_{j,N})]  \\
     =      \frac{1}{P_{N}'(x_{j,N})} \left( \int_\mathbb{R} P_0(x) \, d\mu_2(x)\ [C_1A_1 + C_2 B_1] \right. \\
       \left. + \int_\mathbb{R} P_1(x) \, d\mu_2(x) \  [C_1A_2 + C_2B_2]   \right) . 
\end{multline}
If we denote
\[   D_{1,1} = \int_\mathbb{R} P_0(x) \, d\mu_1(x), \quad   D_{1,2} = \int_\mathbb{R} P_1(x) \, d\mu_1(x), \]
\[   D_{2,1} = \int_\mathbb{R} P_0(x) \, d\mu_2(x), \quad   D_{2,2} = \int_\mathbb{R} P_1(x) \, d\mu_2(x), \]
and use the biorthogonality
\[    \int_\mathbb{R} P_0(x) [A_1 w_1(x) + B_1 w_2(x)]\, d\mu(x) = 1, 
   \quad \int_\mathbb{R} P_0(x) [A_2 w_1(x) + B_2 w_2(x)]\, d\mu(x) = 0, \]
\[  \int_\mathbb{R} P_1(x) [A_1 w_1(x) + B_1 w_2(x)]\, d\mu(x) = 0, 
   \quad \int_\mathbb{R} P_1(x) [A_2 w_1(x) + B_2 w_2(x)]\, d\mu(x) = 1, \]
then we find
\[       \begin{pmatrix} A_1 & 0 \\ A_2 & B_2  \end{pmatrix} \begin{pmatrix} D_{1,1} & D_{1,2} \\ D_{2,1} & D_{2,2} \end{pmatrix} = \begin{pmatrix}  1 & 0 \\ 0 & 1  \end{pmatrix}.  \]
Observe that $A_1=A_{1,0}(x)$, $A_2=A_{1,1}(x)$ and $B_2=B_{1,1}(x)$ are polynomials of degree $0$ and hence these are constants, whereas $B_1 = B_{1,0}=0$ and by the orthogonality of type II multiple orthogonal polynomials we also have $D_{1,2}=0$.

In order to find a formula for $P_N'(x_{j,N})$, we take the limit $x \to x_{j,N}$ in \eqref{CDx-xj} to find
\begin{equation}  \label{CDxjxj} 
\sum_{k=0}^{N-1} P_k(x_{j,N}) [C_1A_{k+1}(x_{j,N})+C_2B_{k+1}(x_{j,N})] 
  = P_{N}'(x_{j,N}) [C_1A_{N}(x_{j,N})+C_2B_{N}(x_{j,N})].   
\end{equation}
The left hand side is the inner product of the left eigenvector and the right eigenvector of $H_N$ for the eigenvalue $x_{j,N}$.
Combining \eqref{lambda1a}, \eqref{lambda2a} and \eqref{CDxjxj} and taking into account
Property \ref{prop5} and \ref{prop6} then gives the main result.
\end{proof}

\section{Examples}    \label{sec6}
\subsection{Bessel functions $K_\nu$}
We will consider the multiple orthogonal polynomials related to the modified Bessel functions $K_{\nu}$ and $K_{\nu+1}$.
Let $(w_1,w_2) = x^\alpha(\rho_\nu,\rho_{\nu+1})$ with
\[    \rho_\nu(x) = 2 x^{\nu/2} K_\nu(2\sqrt{x}), \qquad x >0. \]
We take $\alpha=1$ and $\nu=0$ and we compute the eigenvalues of $H_N$ for $N=10$, together with the weights
$\lambda_{j,10}^{(1)}$ and $\lambda_{j,10}^{(2)}$ using Theorem \ref{thm3}. The type I polynomials
$(A_{1,0},B_{1,0})$ and $(A_{1,1},B_{1,1})$ are the constants
\[   A_{1,0} = 1, \quad B_{1,0}=0, \qquad  A_{1,1} = \alpha+\nu+1, \quad B_{1,1} = -1, \]
and if we use the normalization \eqref{typeInorm}
\[    \int_0^\infty x^\alpha [\widehat{A}_{n,m}(x) \rho_\nu(x) + \widehat{B}_{n,m}(x) \rho_{\nu+1}(x)] x^{n+m-1}\, dx = 1, \]
then we get
\[     \begin{pmatrix} A_1 & 0 \\ A_2 & B_2  \end{pmatrix} = 
\begin{pmatrix} \frac{1}{\Gamma(\alpha+\nu+1)\Gamma(\alpha+1)} & 0 \\
       -\frac{\alpha+\nu+1}{\Gamma(\alpha+\nu+2)\Gamma(\alpha+2)} & \frac{1}{\Gamma(\alpha+\nu+2)\Gamma(\alpha+2)} 
\end{pmatrix}, \]
and hence for $\alpha=1$ and $\nu=0$ we have
\[   \begin{pmatrix} D_{1,1} & 0 \\ D_{2,1} & D_{2,2} \end{pmatrix}
    = \begin{pmatrix}  1 & 0 \\ -1/2 & 1/4  \end{pmatrix}^{-1} 
    = \begin{pmatrix} 1 & 0 \\ 2 & 4  \end{pmatrix}.  \]
We used Maple with \texttt{Digits:=100}
for our calculations and we show the 20 first digits after the decimal point in Fig. \ref{table1}.

\begin{figure}[ht]
\centering
\begin{tabular}{|c|c|c|c|}
\hline
$j$ & $x_{j,10}$ & $\lambda_{j,10}^{(1)}$ &  $\lambda_{j,10}^{(2)}$  \\ 
\hline
  1 &     0.52720348133440875760 & 0.27736269648616286974 & 0.26086734230400106004 \\
  2 &     2.74106066716069179819 & 0.46938499819336417730 & 0.88799214753397210390 \\  
  3 &     8.13937609771412899056 & 0.21135584109286564463 & 0.65379039925659229785 \\  
  4 &   18.66164146312871349710 & 0.03854365644852726770 & 0.17589229666877292663 \\ 
  5 &   36.89653691488348638176 & 0.00322544756122977083 & 0.02038307627872880093 \\  
  6 &   66.43703332978391524587 & 0.00012523808693942895 & 0.00105166051829272396 \\
  7 & 112.55686514754090244347 & 0.00000210903533490802 & 0.00002289663649071884 \\
  8 & 183.67841427499791701294 & 0.00000001307455465436 & 0.00000018043669350953 \\
  9 & 295.27746298319776238423 & 0.00000000002101777610 & 0.00000000036637784733 \\
10 & 485.08440564025807348828 & 0.00000000000000350239 & 0.00000000000007801100 \\
\hline
\end{tabular}
\caption{Quadrature nodes and quadrature weights for $(w_1,w_2)=x(\rho_0,\rho_1)$.}
\label{table1}
\end{figure} 

\medskip
\noindent\textbf{Algorithm 1: Maple}
\begin{Verbatim}[frame=single]
  with(LinearAlgebra):
  b:=n->(n+alpha+1)*(3*n+alpha+2*nu)-(alpha+1)*(nu-1);
  c:=n->n*(n+alpha)*(n+alpha+nu)*(3*n+2*alpha+nu);
  d:=n->n*(n-1)*(n+alpha)*(n+alpha-1)*(n+alpha+nu)*(n+alpha+nu-1);
  H:=n->BandMatrix([[seq(d(k),k=2..n-1)],[seq(c(k),k=1..n-1)],
    [seq(b(k),k=0..n-1)],[seq(1,k=1..n-1)]]);
  alpha:=1.;nu:=0.;
  Digits:=100;
  N:=10;
  E,RE:=Eigenvectors(H(N)):     % right eigenvalues/vectors (descending)
  eig:=Eigenvectors(Transpose(H(N)), output = list):
  eig := sort(eig, (a, b) -> is(Re(b[1]) < Re(a[1]))):    
  EE:=Vector[column]([seq(eig[k][1],k=1..N)]): % eigenvalues (descending)
  LE:=Matrix([seq(op(eig[k][3]), k = 1 .. N)]): % left eigenvectors 
  A:=Matrix([[1/(GAMMA(alpha+nu+1)*GAMMA(alpha+1)),0], 
        [-(alpha+nu+1)/(GAMMA(alpha+nu+2)*GAMMA(alpha+2)),
        1/(GAMMA(alpha+nu+2)*GAMMA(alpha+2))]]);
  D:=MatrixInverse(A);
  u:=j->LE(1:N,j);                     % left eigenvectors
  v:=j->RE(1:N,N-j+1)/RE(1,N-j+1);     % normalized right eigenvectors
  lambda1:=j->D(1,1)*u(j)[1]/DotProduct(u(j),v(j));
  lambda2:=j->(D(2,1)*u(j)[1]+D(2,2)*u(j)[2])/DotProduct(u(j),v(j));
  seq(EE(k),k=1..N);               % quadrature nodes
  seq(lambda1(k),k=1..N);          % quadrature weights (first weight)
  seq(lambda2(k),k=1..N);          % quadrature weights (second weight)
\end{Verbatim}

To illustrate the quality of the quadratures, we will simultaneously compute
\begin{eqnarray*}
   I_1 &=& \int_0^\infty e^{-x} x \rho_0(x)\, dx = 0.1926947246\ldots, \\
   I_2 &=& \int_0^\infty e^{-x} x \rho_1(x)\, dx = 0.2109579130\ldots
\end{eqnarray*}   
using $N=10,20,30,40,50$ nodes in Maple. We show the first 10 decimals in Fig. \ref{tab3}.

\begin{figure}[ht]
\centering
\begin{tabular}{|c|c|c|}
\hline
$N$ &  $I_1$  &  $I_2$  \\
\hline
10 & 0.\textbf{19}40521520      & 0.\textbf{21}14457811 \\
20 & 0.\textbf{1926}653563      & 0.\textbf{2109}395236 \\
30 & 0.\textbf{19269}58911      & 0.\textbf{2109}610461 \\
40 & 0.\textbf{1926947}184      & 0.\textbf{210957}6142  \\
50 & 0.\textbf{1926947}165      & 0.\textbf{21095791}57 \\
\hline
\end{tabular}
\caption{Simultaneous quadrature to $I_1$ and $I_2$}
\label{tab3}
\end{figure}

The results show that the quadrature formulas converge rather slowly to the correct value. A possible  explanation is that
the quadrature nodes become very large and the corresponding quadrature weights are very small, so only
a small proportion of terms in the quadrature sum contribute to the result.

\subsection{Bessel functions $I_\nu$}
For the multiple orthogonal polynomials with weights $(w_1,w_2)=(\omega_{\nu,c},\omega_{\nu+1,c})$, where
\[   \omega_{\nu,c}(x) = x^{\nu/2}I_{\nu}(2\sqrt{x})e^{-cx}, \qquad x > 0,  \]
we take $c=1$ and $\nu=0$ and compute the eigenvalues of $H_N$ for $N=10$,
together with the weights $\lambda_{j,10}^{(1)}$ and $\lambda_{j,10}^{(2)}$ using Theorem \ref{thm3}.
After normalization we get
\[  \begin{pmatrix} A_1 & 0 \\ A_2 & B_2 \end{pmatrix} = e^{-1/c} \begin{pmatrix} c^{\nu+1} & 0 \\ -c^{\nu+2} & c^{\nu+3} \end{pmatrix} , \]
so that for $c=1$ and $\nu=0$ we have
\[    \begin{pmatrix} D_{1,1} & 0 \\ D_{2,1} & D_{2,2} \end{pmatrix} = \begin{pmatrix} 1/e & 0 \\ -1/e & 1/e \end{pmatrix} ^{-1} =
\begin{pmatrix} e & 0 \\ e & e \end{pmatrix}.  \]
We use Matlab and we show the first 10 digits after the decimal point in Fig. \ref{table2}.

\begin{figure}[ht]
\centering
\begin{tabular}{|c|c|c|c|}
\hline
$j$ & $x_{j,10}$ & $\lambda_{j,10}^{(1)}$ &  $\lambda_{j,10}^{(2)}$  \rule{0pt}{14pt} \\
\hline
 1   & 0.1531952228   &  0.3913749988  &  0.0557885974  \\
 2   & 0.8105837014   &  0.8175616919  &  0.4874004644  \\
 3   & 2.0077223654   &  0.8459198767  &  0.9551942639  \\
 4   & 3.7719525634   &  0.4850707607  &  0.8091738873  \\
 5   & 6.1482336073   &  0.1517396396  &  0.3357737316  \\
 6   & 9.2079873838   &  0.0246520172  &  0.0683288497  \\
 7   & 13.0663024491 &  0.0019027391  &  0.0063827530  \\
 8   & 17.9203555594 &  0.0000595495  &  0.0002366956  \\
 9   & 24.1543375116 &  0.0000005543  &  0.0000025816  \\
 10 & 32.7593296369 &  0.0000000007  &  0.0000000038  \\
 \hline
\end{tabular}
\caption{Quadrature nodes and quadrature weights for $(w_1,w_2)=(\omega_{0,1},\omega_{1,1})$.}
\label{table2}
\end{figure}
\newpage

\noindent\textbf{Algorithm 2: Matlab}
\begin{Verbatim}[frame=single]
1   C=1
2   nu=0
3   n=10
4   b=(1+C*(nu+1+2*[0:1:n-1]))/C^2
5   c=[1:1:n-1].*(2+C*(nu+[1:1:n-1]))/C^3
6   d=[2:1:n-1].*([2:1:n-1]-1)/C^4
7   H=diag(ones(1,n-1),1)+diag(b)+diag(c,-1)+diag(d,-2)
8   [R,E,L]=eig(H)      
9   U=L     
10  V=R/diag(R(1,1:n))    
11  D=exp(1/C)*inv([[C^(nu+1),0];[-C^(nu+2),C^(nu+3)]])
12  lambda1=D(1,1)*U(1,1:n)/dot(U,V) 
13  lambda2=(D(2,1)*U(1,1:n)+D(2,2)*U(2,1:n))./dot(U,V)
14  [diag(E),lambda1',lambda2'] 
\end{Verbatim}

We will simultaneously compute the integrals
\begin{eqnarray*}
     J_1 & = & \int_0^\infty \cos(x) \omega_{0,1}(x)\, dx = 0.328224976685277123104160354501976758\ldots \\
     J_2 & = & \int_0^\infty \cos(x) \omega_{1,1}(x)\, dx = -0.39521954160680745592163128352397786234\ldots
\end{eqnarray*}
with $N=10,20,30,40,50$ quadrature nodes in Maple (with \texttt{Digits:=100}). The results are in Fig. \ref{tab4}.
The results in Matlab were comparable for $N=10$ but for $N=20$ and higher some of the eigenvalues and eigenvectors
became complex, due to machine precision. This was avoided in Maple by choosing \texttt{Digits:=100}. Alternatively
one could use variable precision arithmetic (vpa) in Matlab but then Algorithm 2 needs to be modified because
the command \texttt{eig} does not give the left eigenvalues.

\begin{figure}[ht]
\centering
\begin{tabular}{|c|c|c|}
\hline
$N$ &  $J_1$  &  $J_2$  \\
\hline
10 & 0.\textbf{328}340082411357      &   -0.\textbf{395}132567462746 \\
20 & 0.\textbf{32822497}721656944454      & -0.\textbf{3952195}3865314722695 \\
30 & 0.\textbf{32822497668527}696693      & -0.\textbf{39521954160680}6392096 \\
40 & 0.\textbf{32822497668527712310}3734621725   & -0.\textbf{39521954160680745592}554825999940  \\
50 & 0.\textbf{328224976685277123104160354}72   & -0.\textbf{39521954160680745592163128}25809 \\
\hline
\end{tabular}
\caption{Simultaneous quadrature to $J_1$ and $J_2$}
\label{tab4}
\end{figure}

Clearly the results of the quadrature are much better than for the previous example.

\section{Concluding remarks}
We have shown that the nodes and the quadrature weights for simultaneous Gaussian
quadrature can be computed numerically using the eigenvalues and the left and right
eigenvectors of a banded Hessenberg matrix, thus extending the well-known approach
of Golub and Welsch for Gaussian quadrature and symmetric tridiagonal matrices. 
In the present paper we investigated the case of two measures, but the approach works
for any number of measures \cite{JCoussWVA}, however the number of diagonals of the Hessenberg matrix increases. Algorithms for the case of three or more measures
using this Golub-Welsch approach have not been implemented yet. For two and three
measures Tomovi\'c and Stani\'c \cite{Tomovic} published an algorithm in which the quadrature weights are computed by solving some system of equations, and this approach was extended
to more measures by Jovanovi\'c et al. \cite{Jovanovic}. Their matrices are not of Hessenberg
type.

One feature that we noticed is that the computation of the right eigenvectors seems to
be easier and faster than that of the left eigenvectors. The software that we used
(Matlab and Maple) sometimes gives the left eigenvectors in a different order than the
right eigenvectors so that one needs to be careful in pairing up the left and right eigenvectors.
For nonsymmetric tridiagonal matrices Van Dooren et al. \cite{VanDooren} gave some improvements
for computing the left and right eigenvectors, and a similar idea might also be valuable
for Hessenberg matrices. These authors already worked on the Hessenberg case in \cite{Laudadio}
and noticed that the total positivity of the matrices allows a stable and efficient way to compute
the eigenvalues. The Hessenberg matrices in our examples \S 6.1 and \S 6.2 are indeed totally
positive as they can be factored using bidiagonal matrices with positive entries, see Lima \cite{Lima}.

The Bessel weights in our example have the remarkable property that one knows the 
recurrence coefficients of the multiple orthogonal polynomials on the stepline exactly with simple
formulas. This is in sharp contrast with the recurrence coefficients of the corresponding
orthogonal polynomials, for which no expression is known, see \cite{Gautschi} for
weight functions related to the  $K_\nu$ Bessel function. This means that we don't need an 
additional step to compute the recurrence coefficients from the moments or modified moments
of the measures, which is usually a badly conditioned problem. We can start right away from
the known recurrence coefficients to build the Hessenberg matrix. In case one knows the
nearest neighbor recurrence coefficients of the multiple orthogonal polynomials, one can compute
the recurrence coefficients for the polynomials on the stepline using an algorithm of Filipuk et al.
\cite{FilHanWVA}.

\section{Declarations}
\subsection{Ethical approval}
Not applicable.
\subsection{Availability of supporting data}
All the data are numerical results published in the paper. The algorithms are also in the paper.
\subsection{Competing interests}
Not applicable.
\subsection{Funding}
Supported by research project G0C9819N of FWO (Research Foundation -- Flanders).
\subsection{Author's contribution}
This is a single author paper. The paper was written by the author and the numerical results were obtained by the author.
\subsection{Acknowledgments}
The author is grateful to the referees for useful comments and additional references.


\begin{verbatim}
Walter Van Assche
Department of Mathematics
KU Leuven
Celestijnenlaan 200B, box 2400
BE-3001 Leuven, BELGIUM
walter.vanassche@kuleuven.be
orcid 0000-0003-3446-6936
\end{verbatim}

\end{document}